\documentclass[hidelinks]{article}
\usepackage{graphicx} 
\usepackage{hyperref}
\usepackage{amsthm}
\usepackage{thm-restate}
\usepackage{geometry}
\usepackage{amssymb}
\usepackage{amsmath}
\usepackage{cleveref}
\usepackage{comment}
\usepackage{xcolor}
\usepackage[shortlabels]{enumitem}
\usepackage{listings}
\usepackage{blkarray}
\usepackage{lmodern}
\usepackage[english]{babel}
\usepackage{float}
\usepackage{tikz}

\newtheorem{theorem}{Theorem}[section]
\newtheorem{lemma}[theorem]{Lemma}

\newtheorem{claim}[theorem]{Claim}

\newtheorem{conjecture}[theorem]{Conjecture}

\theoremstyle{definition}

\newcommand{\cI}{\mathcal{I}}

\newcommand{\cH}{\mathcal{H}}

\newcommand{\bracenom}{\genfrac{\lbrace}{\rbrace}{0pt}{}}

\title{Hypergraph universality via branching random walks}

\author{
Rajko Nenadov\thanks{School of Computer Science, University of Auckland, New Zealand. Email: \texttt{rajko.nenadov@auckland.ac.nz}. Research supported by the Marsden Fund of the Royal Society of New Zealand.}
}
\date{}

\begin{document}

\maketitle

\begin{abstract}
Given a family of hypergraphs $\mathcal{H}$, we say that a hypergraph $\Gamma$ is $\mathcal{H}$-universal if it contains every $H \in \mathcal{H}$ as a subgraph. For $D, r \in \mathbb{N}$, we construct an $r$-uniform hypergraph with $\Theta\left(n^{r - r/D} \log^{r/D}(n)\right)$ edges which is universal for the family of all $r$-uniform hypergraphs with $n$ vertices and maximum degree at most $D$. This almost matches a trivial lower bound $\Omega(n^{r - r/D})$ coming from the number of such hypergraphs. 

On a high level, we follow the strategy of Alon and Capalbo used in the graph case, that is $r = 2$. The construction of $\Gamma$ is deterministic and based on a bespoke product of expanders, whereas showing that $\Gamma$ is universal is probabilistic. Two key new ingredients are a decomposition result for hypergraphs of bounded density, based on Edmond's matroid partitioning theorem, and a tail bound for branching random walks on expanders.
\end{abstract}

\section{Introduction}

A graph $\Gamma$ is universal for a family of graphs $\mathcal{H}$ if every $H \in \mathcal{H}$ is a subgraph of $\Gamma$, not necessarily induced. Constructions of sparse graphs which are universal for families $\mathcal{H}$ of interest has been studied extensively in the past decades. This includes families of graphs with bounded degree \cite{alon02sparse,alon07sparsebounded,alon08optimal,alon2000universality,alon2001near}, which are also the central family considered in this paper, trees and forests \cite{chung78tree,chung78tree2,chung83spanning,friedman87tree} and, more generally, graphs with bounded degeneracy~\cite{allen2023universality,nenadov2016ramsey}, as well as families of graphs with additional structural properties such as planar graphs \cite{babai82planar,esperet23planar} and, more generally, graphs with small separators \cite{capalbo02planar,capalbo99small,chung90separator}. Other than being an interesting problem in its own right, motivation for studying sparse universal graphs comes from applications in VLSI circuit design \cite{bhatt84vlsi}, data storage \cite{chung83storage}, and simulation of parallel computer architecture \cite{bhatt1986optimal}, to name a few.

We are interested in the family $\mathcal{H}^{(r)}(D, n)$ consisting of all $r$-uniform hypergraphs with $n$ vertices and maximum degree at most $D$. The case $r = 2$ was first studied by Alon, Capalbo, Kohayakawa, R\"odl, Ruci\'nski, and Szemer\'edi \cite{alon2000universality}, where they constructed an $\mathcal{H}^{(2)}(D, n)$-universal graph with $\Theta\left(n^{2 - 1/D} \log^{1/D}(n) \right)$ edges. This was improved in a series of papers \cite{alon02sparse,alon07sparsebounded,alon2001near}, culminating with the work of Alon and Capalbo \cite{alon08optimal} where they constructed a universal graph with $\Theta\left( n^{2 - 2/D} \right)$ edges. A simple counting argument based on the size of the family $\mathcal{H}^{(2)}(D, n)$ gives a lower bound $\Omega \left( n^{2 - 2/D} \right)$, showing that the construction of Alon and Capalbo is optimal.

The hypergraph case $r \ge 3$ was considered by Parczyk and Person \cite{parczyk16hypergraphs} and Hetterich, Parczyk, and Person \cite{hetterich16hypergraphs}. By reducing the problem to the graph case, they showed that for even $r$ there exists an $r$-uniform hypergraph ($r$-graph for short) with $\Theta\left( n^{r - r/D} \right)$ edges which is $\mathcal{H}^{(r)}(D, n)$-universal. They also obtained the same bound in the case $r$ is odd and $D = 2$. Similarly as in the graph case, a simple counting argument based on the size of the family $\mathcal{H}^{(r)}(D, n)$ shows that this is optimal. In the case of odd $r$ and $D > 2$, they constructed a universal hypergraph with $\Theta\left( n^{r - (r+1)/\Delta} \right)$ edges, where $\Delta = \lceil (r+1)D / r \rceil$. This falls short of the lower bound $\Omega\left( n^{r - r/D} \right)$ by a polynomial factor for every such $r$ and $D$. Our main result is the existence of universal $r$-graphs which are off only by a small logarithmic factor.

\begin{theorem} \label{thm:bounded}
For every $r, D \in \mathbb{N}$ there exists $C = C(r, D) > 0$ such that the following holds: For every $n \in \mathbb{N}$ there exists an $r$-graph $\Gamma$ with
$$
    e(\Gamma) \le C n^{r - r / D} \log^{r/D}(n),
$$
which is $\mathcal{H}^{(r)}(D, n)$-universal.
\end{theorem}

On a high-level, our proof follows the construction based on a product of expander graphs used by Alon and Capalbo \cite{alon07sparsebounded,alon08optimal}. To show that the obtained hypergraph is $\mathcal{H}^{(r)}(D, n)$-universal, we employ two new ingredients. The first one, Lemma \ref{lemma:hyp_decomposition}, provides a collection of graphs with a simple structure (namely unicyclic) which together \emph{underpin} a given hypergraph $H$. The construction of $\Gamma$ is tailored to make use of this. The second one is a first step towards generalising a classic result of Gilman \cite{gilam98chernoff} on tail bounds for random walks on expanders to branching random walks which follow a given `blueprint' tree with bounded degree. This is used to show that a combination of certain graph homomorphisms together form an injection -- and thus an embedding of $H$ in $\Gamma$. We consider this tail bound for branching random walks to be our main technical contribution, and an interesting research direction in its own right, thus we describe it next. 

\subsection{Tail bound for branching random walks}
\label{sec:branching_walk_intro}

In the simplest form, a random walk of length $\ell$ on a graph $G$ is defined as follows: Let $v_1$ be a vertex in $V(G)$ chosen uniformly at random, and for each $i \in \{2, \ldots, \ell-1\}$, sequentially, take $v_i$ to be a neighbour of $v_{i-1}$ chosen uniformly at random. In the case $G$ has bounded degree, this gives an efficient way of sampling $\ell$ vertices from $G$ in terms of the number of random bits. Indeed, compared to $\ell \log(|V(G)|)$ bits required to sample $\ell$ vertices completely uniformly and independently, random walk requires only $\log(|V(G)|) + (\ell - 1) \log(D)$, where $D$ is the maximum degree in $G$. This is of great importance in theoretical computer science, and it prompts the question of how much sampling vertices using a random walk resembles the uniform distribution. While the vertices in a random walk are very much dependent locally, it turns out that globally the two distributions exhibit similar phenomena. This was first observed by Ajtai, Koml\'os, and Szemer\'edi \cite{aks87} who studied the probability that a random walk stays confined to a given subset. Their result was significantly strengthened by Gillman \cite{gilam98chernoff}, who showed that that if $G$ is a good expander, then the probability  that a random walk hits a given set $S \subseteq V(G)$ significantly more that $\ell |S| / |V(G)|$ times is similar to what Chernoff-Hoeffding inequality gives when all vertices are sampled uniformly and independently. Since then, there has been a great interest in generalising these results and by now there is a large body of research on tail bounds for random walks on expander graphs and, more generally, finite state Markov chains (e.g.~see \cite{dinwoodie95,garg18matrix,hea08sampling,kahale97,leon04,lez98,naor20boundedmoment,rao17sharp,wagner08}).

We propose the study of similar questions for a certain class of branching random walks. For a rooted tree $T$ and a graph $G$, we define a \emph{random $T$-walk} on $G$ to be a homomorphism $\phi \colon T \to G$ given by the following random process. Let $t_1, \ldots, t_N$ be any ordering of the vertices in $T$ such that $t_1$ is the root, and if $t_i$ is closer to the root than $t_j$ then $i < j$. Choose $\phi(t_1)$ uniformly at random from $V(G)$, and for $i \ge 2$, sequentially, choose $\phi(t_i)$ uniformly at random among the set of neighbours of $\phi(p_i)$, where $p_i$ is the parent of $t_i$ (thus precedes it in the considered ordering). Note that the choice of the ordering is irrelevant for the outcome of the process, as long as it satisfies the stated property. In this terminology, a random walk of length $\ell-1$ corresponds to a random $P_\ell$-walk where $P_\ell$ is a path with $\ell$ vertices.  

Recall that an $(n,d,\lambda)$-graph is a $d$-regular graph with $n$ vertices, such that the second largest absolute value of its adjacency matrix is at most $\lambda$. The following lemma is our main technical contribution.

\begin{lemma} \label{lemma:random_branching}
    There exist constants $K_0 > 1$ and  $\alpha > 0$ such that the following holds. Let $T$ be a rooted tree with maximum degree $D$, and suppose $G$ is an $(n, d, \lambda)$-graph with $\lambda  < \alpha d / D$. Given a non-empty subset $U \subseteq V(T)$ and a vertex $x \in V(G)$, let $X$ denote the number of vertices from $U$ which are mapped to $x$ in a random $T$-walk on $G$. Then $\mathbb{E}[X] = |U| / n$ and
    $$
        \Pr\left[X > K \mathbb{E}[X] \right] \le e^{-(K-K_0) \mathbb{E}[X]},
    $$
    for every $K > K_0$.
\end{lemma}

On the one hand, Lemma \ref{lemma:random_branching} is weaker than the previously discussed result of Gillman \cite{gilam98chernoff} in two ways: (i) we require $K$ to be sufficiently large, and (ii) we only bound the number of times a particular vertex $x \in V(G)$ is being hit, rather than a subset of vertices.  On the other hand, it is stronger in the sense that it supports any tree $T$ and not just a path, and it also counts the number of times a particular subset of vertices of $T$ hits $x$, instead of the whole of $T$. In the case of random walks, this corresponds to the number of times a particular set of steps has landed on $x$. A common generalisation would be to consider a set of functions $f_{t}: V(G) \to [0,1]$, one for each pair $t \in V(T)$, and establish a tail bound for the random variable $X = \sum_{t \in V(T)} f_t(\phi(t))$. In the case of random walks, this has been done by Rao and Regev \cite{rao17sharp}. Another, rather obvious, open problem is to improve the lower bound on $K_0$ for which the conclusion of Lemma \ref{lemma:random_branching} holds. We leave these as interesting directions for future research.

\section{Decomposition Lemma} \label{sec:decomposition}

Given a hypergraph $H$, recall the usual definition of the maximal \emph{density} of $H$:
$$
    m(H) = \max_{H' \subseteq H} \frac{e(H')}{v(H')}.
$$

\begin{lemma} \label{lemma:hyp_decomposition}
Let $H$ be an $r$-graph with $m(H) \le a / b$, for some $a, b \in \mathbb{N}$ with $b > r-1$ and $a > b / (r-1)$. Then there exists a family of graphs (that is, $2$-graphs) $H_1, \ldots, H_a$ on the vertex set $V(H)$ such that the following holds:
\begin{enumerate} [(D1)]
    \item \label{prop:unic} Each connected component of every $H_i$ is unicyclic and the maximum degree of $H_i$ is at most $2D$, where $D$ is the maximum degree in $H$, and 
    \item \label{prop:forest} For each hyperedge $h \in E(H)$ there exist forests $F_1^{(h)} \subseteq H_1, \ldots, F_a^{(h)} \subseteq H_a$ on the vertex set $h$ (recall that $h \subseteq V(H)$ is a subset of size $r$) such that $\sum_{i = 1}^a e(F_i^{(h)}) = b$.
\end{enumerate}
\end{lemma}

It is crucial that each $F_i^{(h)}$ is a forest and not just a unicyclic graph, like what we require in \ref{prop:unic}. The proof of Lemma \ref{lemma:hyp_decomposition} uses Edmonds' matroid partitioning theorem \cite{edmonds65matroid}. Recall that a finite matroid $M$ is a pair $(E, \mathcal{I})$ where $E$ is a finite set and $\mathcal{I}$ is a family of subsets of $E$ with the following properties:
\begin{itemize}
    \item $\emptyset \in \mathcal{I}$,
    \item If $A' \subseteq A \subseteq E$ and $A \in \mathcal{I}$, then $A' \in \mathcal{I}$, and
    \item If $A, B \in \mathcal{I}$ and $|A| > |B|$, then there exists $x \in A \setminus B$ such that $B \cup \{x\} \in \mathcal{I}$.
\end{itemize}
The set in $\mathcal{I}$ are referred to as \emph{independent sets}.

\begin{theorem}[Edmonds' partitioning theorem]
Let $M = (E, \mathcal{I})$ be a finite matroid, and let
$$
    k(M) = \max_{S \subseteq E} \left\lceil \frac{|S|}{r(S)} \right\rceil,
$$
where $r(S)$ denotes the size of a largest independent set from $\mathcal{I}$ which is contained in $S$. Then there exists a partition $E = I_1 \cup I_2 \cup \ldots \cup I_{k(M)}$ such that $I_i \in \mathcal{I}$ for each $i \in \{1, \ldots, k(M)\}$. 
\end{theorem}

We are now ready to prove Lemma \ref{lemma:hyp_decomposition}.

\begin{proof}[Proof of Lemma \ref{lemma:hyp_decomposition}]
Given an $r$-uniform $H$ with $m(H) \le a / b$, we construct a bipartite graph $B$ on vertex sets $U$ and $V$ as follows. The set $V$ corresponds to the vertex set $V(H)$, and for each hyperedge $h \in H$ there are $b$ vertices in $U$ corresponding to $h$. We put an edge between $u \in U$ and $v \in V$ iff the hyperedge in $H$ corresponding to $u$ contains $v$. 

Let $\cI$ be the family of all subsets $X \subseteq U$ such that: (i) $X$ contains at most $r-1$ vertices corresponding to each hyperedge, and (ii) for each nonempty $X' \subseteq X$ we have $|N_B(X')| \ge |X'|$. We claim that $M = (U, \cI)$ is a matroid with $\cI$ being the family of independent sets. We trivially have $\emptyset \in \cI$, and by the definition we have that $Y \in \cI$ and $X \subseteq Y$ implies $X \in \cI$. It remains to verify the \emph{augmentation} axiom, that is, for each $X, Y \in \cI$ such that $|X| < |Y|$, there exists $x \in Y \setminus X$ such that $X \cup \{x\} \in \cI$. This can be seen as follows. We can assume that if $X_h \subseteq X$ and $Y_h \subseteq Y$ denote subsets which corresponds to vertices associated with a hyperedge $h$ and $|X_h| \le |Y_h|$, then $X_h \subseteq Y_h$. This is because each vertex in $Y_h$ has the same neighbourhood, thus we can reassign vertices such that this holds. By the definition of $\cI$ and Hall's condition, there exist matchings $M_X$ and $M_Y$ in $B$ which saturate $X$ and $Y$, respectively. Now $M_X \cup M_Y$ contains an augmenting path which gives a matching saturating $X \cup \{x\}$ for some $x \in Y \setminus X$. Therefore, $|N(X')| \ge |X'|$ for every $X' \subseteq X \cup \{x\}$. By the initial assumption, we have that the number of vertices in $X \cup \{x\}$ corresponding to each hyperedge is, again, at most $r-1$, thus $X \cup \{x\} \in \cI$.

Now that we have established that $M$ is a matroid, we can state a claim that is the heart of the proof of the lemma.

\begin{claim} \label{claim:Mf_partition}
    There exists a partition $U = U_1 \cup \ldots \cup U_a$ such that $U_1, \ldots, U_a \in \cI$.
\end{claim}
\begin{proof}[Proof of Claim \ref{claim:Mf_partition}]
    By Edmonds' partition theorem, it suffices to show $|Z| / r(Z) \le a$ for every $Z \subseteq U$, where $r(Z)$ denotes the \emph{rank} of $Z$, that is, the size of a largest independent set from $M$ which is contained in $Z$.

    Consider some $Z \subseteq U$, and let $Z' \subseteq Z$ be obtained from $Z$ by removing all but at most $r-1$ vertices corresponding to each hyperedge from $H$. Then $r(Z)$ is equal to the size of a largest matching in $B$ between $Z'$ and $V$. By K\"{o}nig's theorem, $r(Z) = |C|$ where $C$ is a smallest vertex cover of the induced bipartite subgraph $B[Z', V]$. Let $Z_C = Z' \cap C$, $\hat Z_C = Z' \setminus C$, 
    $V_C = V \cap C$, and $\hat V_C = V \setminus V_C$. Note that there is no edge in $B$ between $\hat Z_C$ and $\hat V_C$. Moreover, the smallest vertex cover in $B[\hat Z_C, V_C]$ is of size $|V_C|$, and in $B[Z_C, \hat V_C]$ of size $|Z_C|$. Therefore, a largest matching between $\hat Z_C$ and $V_C$ is of size $|V_C|$, and a largest matching between between $Z_C$ and $\hat V_C$ is of size $|Z_C|$. As every hyperedge corresponding to a vertex in $\hat Z_C$ is fully contained in $V_C$, we conclude
    \begin{equation} \label{eq:hat_ZC}
        |\hat Z_C| \le (r-1) \cdot m(H) |V_C| \le \frac{(r-1)a}{b} |V_C|.
    \end{equation}
    We have $|Z| \le b |Z'| / (r-1)$, thus
    $$
        \frac{|Z|}{r(Z)} = \frac{|Z|}{|Z_C| + |V_C|} \le \frac{b}{r-1} \cdot \frac{|Z'|}{|Z_C| + |V_C|} = \frac{b}{r-1} \cdot \frac{|Z_C| + |\hat Z_C|}{|Z_C| + |V_C|} \le \frac{b}{r-1} \cdot \max\{ 1, |\hat Z_C| / |V_C| \}.
    $$
    From $b / (r-1) < a$ and \eqref{eq:hat_ZC}, we conclude $|Z| / r(Z) \le a$, as desired.
\end{proof}

For each hyperedge $h \in H$ fix a cyclic ordering of its vertices. Let us denote with $h(v)$ the successor of a vertex $v \in h$ in such an ordering for a hyperedge $h$. By the definition of $\cI$ and Hall's condition, for every independent set $X \in \cI$ there exists a matching in $B$ which saturates $X$. Let $\phi \colon U_i \rightarrow V$ denote such a matching saturating $U_i$. We form $H_i$ by taking an edge $\{\phi(u), h_u(\phi(u))\}$ for each $u \in U_i$, where $h_u \in H$ is the hyperedge in $H$ corresponding to $u$. A vertex is incident to at most two edges coming from each hyperedge it is part of, and every connected component of $H_i$ contains at most one cycle (one can think of the obtained graph as being a directed graph with out-degree at most $1$), thus \ref{prop:unic} is satisifed.A forest $F_i^{(h)}$ corresponding to the hyperedge $h \in H$ is simply a (possibly empty) collection of paths given by the union of edges 
$\{\phi(u), h(\phi(u))\}$ for $u \in U_i$ corresponding to $h$. Note that this is indeed a forest, and not a cycle, as $U_i$ contains at most $r-1$ vertices from $U$ corresponding to $h$. Each vertex in $U$ corresponding to $h$ contributes exactly one edge to some $F_i^{(h)}$, thus \ref{prop:forest} holds as well.
\end{proof}

\section{Branching random walk on expanders}

In this section we prove Lemma \ref{lemma:random_branching}. The proof follows the strategy of Rao and Regev \cite{rao17sharp}, with the following lemma being the key new ingredient. This lemma is also the main difference compared to \cite{rao17sharp}, which deals with the simpler case of random $P_\ell$-walks. 
\begin{lemma} \label{lemma:key}
    Let $T$ be a rooted tree with maximum degree $D$, and suppose $G$ is an $(n, d, \lambda)$-graph, for some $\lambda < d / (2^{10} D)$. Consider a random $T$-walk on $G$. For a subset $W \subseteq V(T)$ and $x \in V(G)$, let $I_x(W)$ be the indicator random variable for the event that all the vertices in $W$ are mapped to $x$. Then for any $x \in V(G)$, a subset $U \subseteq V(T)$, and $1 \le k \le |U|$, we have
    \begin{equation} \label{eq:sum_k}
        \mathbb{E}\left[ \sum_{W \in \binom{U}{k}} I_x(W) \right] \le \sum_{i = 1}^{k} \binom{k-1}{i-1} \frac{(2^8|U|/n)^{i}}{i!} \left( 2^9 D \; \frac{\lambda}{d} \right)^{k-i},
    \end{equation}
    where $\binom{U}{k}$ denotes the family of all $k$-element subsets of $U$.
\end{lemma}

In the proof of Lemma \ref{lemma:key} we use the following well known property of random walks on expanders, see, e.g., \cite{hoory06expander}.

\begin{lemma} \label{lemma:random_walk_expander} Let $G$ be an $(n, d, \lambda)$-graph, and consider a random walk starting in a given vertex $v \in V(G)$. The probability that after exactly $\ell$ steps we finish in a vertex $w \in V(G)$ is at most
$$
    1/n + (\lambda/d)^\ell.
$$
\end{lemma}

The proof of Lemma \ref{lemma:key} is combinatorial in nature, based on a careful encoding of a depth-first search traversal of the given tree.

\begin{proof}[Proof of Lemma \ref{lemma:key}]
Let us denote with $r$ the root of $T$. For a subset $X \subseteq V(T)$, let $X^\uparrow \subseteq V(T)$ denote the set of all vertices $v \in T$ which lie on the path from some $v \in X$ to the root $r$, including $X$ and $r$ (that is, all vertices `above' $X$ in a top-down drawing of $T$).

Let us first describe the overall strategy of the proof and establish some important notation. For each $W \in \binom{U}{k}$, we define the ordering $\sigma(W) = (w_0, \ldots, w_{k-1})$ of $W$ such that $w_i \not \in W_i^\uparrow$ for every $i \in \{1, \ldots, k-1\}$, where $W_{i} = \{w_0, \ldots, w_{i-1}\}$. For example, taking $\sigma(W)$ to be an ordering induced by the distance of a vertex from the root, tie-breaking in some arbitrary way, satisfies this property. However, it will be important for us that $\sigma(W)$ can be encoded efficiently, thus our algorithm for producing $\sigma$ is more involved and based on depth-first search. We postpone this until it becomes relevant. Next, let $h_i \in W_i^\uparrow$ be the closest vertex in $W_i^\uparrow$ to $w_i$, and let $d_i$ denote their distance. Note that $d_i$ is a function of $i$ and $W$ (that is, of $\sigma(W)$), thus we write $d_i(W)$ to signify this. Conditioned on the outcome of the $T[W_i^{\uparrow}]$-walk $\phi$, $w_i$ is mapped onto the last vertex in a random walk of length $d_i$ starting from $\phi(h_i)$. Note that $\Pr[\phi(w_0) = x] = 1/n$ as the stationary distribution is uniform, thus by Lemma \ref{lemma:random_walk_expander} we have
\begin{equation} \label{eq:expectation_W}
    \mathbb{E}[I_x(W)] = \Pr[I_x(W) = 1] \le \frac{1}{n} \prod_{i = 1}^{k-1} \left(\frac{1}{n} + (\lambda / d)^{d_i} \right).
\end{equation}
We now use the trick of Rao and Regev \cite{rao17sharp} and `unroll' the right hand side,
\begin{equation} \label{eq:estimate_IxW}
    \mathbb{E}[I_x(W)] \le \sum_{f \in \{0,1\}^{k-1}} \left( \frac{1}{n} \right)^{k - |f|} \prod_{i\colon f_i = 1} (\lambda / d)^{d_i},
\end{equation}
where $|f|$ denotes the number of 1's in $f = (f_1, \ldots, f_{k-1})$. Consider some $f \in \{0, 1\}^{k-1}$, and let $W_{\bar{f}} = \{w_0\} \cup \{w_i \colon f_i = 0\}$. Now fix some $F \in \binom{U}{k - |f|}$, and iterate over all sets $W \in \binom{U}{k}$ such that $W_{\bar{f}} = F$. Our goal is to show
\begin{equation} \label{eq:estimate_f}
    \sum_{\substack{W \in \binom{U}{k} \\ W_{\bar{f}} = F}} \prod_{i \colon f_i = 1} (\lambda/d)^{d_i(W)} \le 256^{k - |f|} \left( 2^9 D \; \frac{\lambda}{d} \right)^{|f|}.
\end{equation}
Together with \eqref{eq:estimate_IxW}, this implies the desired bound:
\begin{align*}
    \mathbb{E}\left[ \sum_{W \in \binom{U}{k}} I_x(W) \right] &\stackrel{\eqref{eq:estimate_IxW}}{\le} \sum_{W \in \binom{U}{k}} \sum_{f \in \{0,1\}^{k-1}} \left( \frac{1}{n} \right)^{k - |f|} \prod_{i \colon f_i = 1} (\lambda/d)^{d_i(W)} \\
    &= \sum_{f \in \{0,1\}^{k-1}} \left(\frac{1}{n} \right)^{k - |f|} \sum_{F \in \binom{U}{k-|f|}} \sum_{\substack{W \in \binom{U}{k} \\ W_{\bar{f}} = F}} \prod_{i \colon f_i = 1} (\lambda/d)^{d_i(W)} \\
    &\stackrel{\eqref{eq:estimate_f}}{\leq} \sum_{f \in \{0,1\}^{k-1}} \left(\frac{2^8}{n} \right)^{k - |f|} \binom{|U|}{k - |f|} \left( 2^9 D \; \frac{\lambda}{d}  \right)^{|f|} \\
    &\leq \sum_{i = 0}^{k-1} \binom{k-1}{i} \frac{(2^8 |U| / n)^{k-i}}{(k-i)!} \left( 2^9 D \; \frac{\lambda}{d} \right)^{i}.
\end{align*}
After changing the indexing, we obtain \eqref{eq:sum_k}.

Let us now outline the strategy for estimating \eqref{eq:estimate_f}. For a fixed $f \in \{0, 1\}^{k-1}$ with $|f| \ge 1$, let $I_f = \{i \colon f_i = 1\}$. Our aim is to show that there is an injection $\psi_{f,F}$ from the family of all $W \in \binom{U}{k}$ such that $W_{\bar{f}} = F$ to the set of $4$-tuples $(B, R, (b_i)_{i \in I_f}, (\mathbf{c}_i)_{i \in I_f})$, where $B, R \subseteq \{1, \ldots, 4k\}$, $b_i \in \mathbb{N}_0$, and $\mathbf{c}_i$ is a sequence of finite length with each element being in $\{1, \ldots, D\}$. The mapping $\psi_{f,F}$ also has the following important property: if $\psi_{f,F}(W) = (B, R, (b_i)_{i \in I_f}, (\mathbf{c}_i)_{i \in I_f})$, then $d_i(W) = |\mathbf{c}_i| \ge \max(1, b_i)$ for every $i \in I_f$. We are not yet in a position to say where such a $4$-tuple comes from, other than hint the reader that it is an encoding of a certain traversal of the vertices in $W$ within $T$. Assuming we have such an injection $\psi_{f,F}$, we easily obtain \eqref{eq:estimate_f}:
\begin{align*}
    \sum_{\substack{W \in \binom{U}{k} \\ W_{f} = F}} \prod_{i \colon f_i = 1} (\lambda/d)^{d_i(W)} &\le 4^{4k}\left( \sum_{b_i = 1}^{\infty} \sum_{d_i = b_i}^{\infty} D^{d_i} (\lambda / d)^{d_i}  \right)^{|f|} \\
    &=4^{4k} \left( \sum_{b_i = 1}^\infty \frac{\kappa^{b_i}}{1 - \kappa} \right)^{|f|}
    \le 256^{k - |f|} \left( \frac{256 \kappa}{1 - 2\kappa} \right)^{|f|},
\end{align*}
where $\kappa = D \lambda / d < 1/4$. Now \eqref{eq:estimate_f} follows from $1 - 2\kappa > 1/2$. The importance of the property $b_i \le |\mathbf{c}_i| = d_i(W)$ is evident from the calculation.

It remains to describe $\psi_{f,F}$. To do so, we first describe the algorithm which produces the ordering $\sigma(W)$.

\paragraph{The ordering $\sigma(W)$.} Fix an arbitrary ordering $\pi$ of the vertex set of $T$, and for each vertex $v \in T$ fix an ordering $\pi_v$ on the set of children of $v$. Whenever we refer to the $k$-th child of a vertex $v \in T$, we mean $k$-th according to $\pi_v$. We denote with $T_v$ the subtree of $T$ rooted in $v$, which we identify with its set of vertices. 

Consider some $W \in \binom{U}{k}$. We define the ordering $\sigma(W) = (w_0, \ldots, w_{k-1})$ on $W$ using the depth-first search over a portion of $T$, with a very specific choice of the next vertex to visit. Throughout the procedure we maintain a number of sets, sequences, and indices, initially set as $S = (r)$, $\widehat W = W$, $B = R = \emptyset$, and $j = 0$. If $r \in W$, then set $w_0 = r$ and $j = 1$, and remove $r$ from $\widehat W$. Intuitively, $\widehat W$ is the set of vertices in $W$ which we have not yet visited, $S$ is the stack which keeps track of important vertices on the path from the root to the current vertex, $B$ tells us whether in a particular step we continued the exploration in the subtree ``below'' the current vertex, and $R$ tells us whether in a particular step in which we removed a vertex from the stack, we stayed in the subtree of the next vertex from the top of the stack (i.e.\ we moved to the ``right'' of the current vertex). Finally, whenever we move to a new vertex $w_j$, we record how we got to it in terms of a number $b_i$ -- recording how many steps to go ``back'' towards the root from the current vertex -- and a sequence $\mathbf{c}_i$ -- recording the sequence of moves which tell us how to move through the part of $T$ which has not yet been explored, in order to reach $w_i$.  

Throughout the procedure we use $t$ to denote the ordinal number of the current iteration. As long as $S$ is not empty, repeat:
\begin{enumerate}[(1)]
    \item Let $s$ be the last vertex in $S$ (that is, the top of the stack).
    
    \item If $T_s \cap \widehat W \neq \emptyset$, add $t$ to the set $B$ and choose $w_j \in T_s \cap \widehat W$ to be the closest vertex to $s$, tie-breaking according to $\pi$. Let $d_j$ denote the distance from $s$ to $w_j$ and $\mathbf{c}_j \in [D]^{d_j}$ a sequence describing how to get from $s$ to $u$ in $T$ (recall that $T$ has maximum degree at most $D$ and the children of each node are ordered). For completeness, set $b_j = 0$. Add $w_j$ to the end of $S$, and remove it from $\widehat W$. Increase $j$, and proceed to the next round.
    
    \item Otherwise, we have $T_s \cap \widehat W = \emptyset$. 
    \begin{itemize}
        \item Remove $s$ from $S$. If $S$ is now empty terminate the procedure. 
        \item Let $s'$ be the new vertex on the top of $S$ (that is, the end of $S$). If $T_{s'} \cap \widehat W = \emptyset$, proceed to the next round.
        \item Otherwise, add $t$ to $R$. For each $u \in T_{s'} \cap \widehat W$, let $h(u)$ denote the vertex on the path from $s'$ to $s$ which is closest to $u$. Note that it cannot happen that $h(u) = s$, though it could be $h(u) = s'$. 
        \item Take $u \in T_{s'} \cap \widehat W$ with $h(u)$ closest to $s$ (i.e.\ furthest away from $s'$). If there are multiple such vertices take the one which itself is closest to $s$, that is, the one for which the path from $h(u)$ to $u$ is shortest (tie breaking according to $\pi$). 
        \item If $h(u) \neq s'$ then add $h(u)$ to the end of $S$. 
        \item Add $u$ to the end of $S$ and set $w_j = u$. Set $b_j$ to be the distance from $s$ to $h(u)$ (denoting how many steps we need to go ``back'' towards $s'$), $d_j$ to be the distance from $h(u)$ to $u$, and $\mathbf{c}_j \in [D]^{d_j}$ the description of how to get from $h(u)$ to $u$ in $T$. Increase $j$, and remove $u$ from $\widehat W$. Proceed to the next round.
        \end{itemize}
\end{enumerate}

Observe the following two crucial properties. First, for every $j \in \{0, \ldots, k-1\}$ for which $w_j$ was defined in (3) we have $b_j \le d_j$. If this was not the case, then we would have chosen $w_j$ before $w_{j-1}$. Second, the procedure terminates after $t \le 4k$ rounds.

\paragraph{The mapping $\psi_{f,F}$.} Fix $f \in \{0, 1\}^{k-1}$ with $|f| \ge 1$ and $F \subseteq \binom{U}{k - |f|}$. Consider some $W \in \binom{U}{k}$ such that $W_{\bar{f}} = F$. Run the previously described algorithm on $W$, and let $B, R \subseteq \{1, \ldots, 4k\}$, $(b_i)_{i \in [k]}$, and $(\mathbf{c}_i)_{i \in [k]}$ be as given by the algorithm upon its termination. Define 
$$
    \psi_{f,F}(W) := (B, R, (b_i)_{i \in I_f}, (\mathbf{c}_i)_{i \in I_f}).
$$
The intuition here is that if we know $f$ and $F = W_{\bar{f}}$, and for each $i \in \{1, \ldots, k-1\}$ such that $f_i = 1$ we have enough information on how to reach $w_i$, which is precisely what is encoded in $\psi_{f,F}(W)$, then we can uniquely reconstruct the whole set $W$. Therefore, $\psi_{f,F}$ is necessarily injective. 

To reconstruct $W$ from $f$, $F = W_{\bar{f}}$, and $\psi_{f,F}(W)$, we repeat the algorithm for determining $\sigma(W)$ with the following modification: in the steps corresponding to choosing $w_i$ for $i \in I_f$, instead of taking a vertex from $W$ we follow steps described by $b_i$ and $\mathbf{c}_i$. The sets $B$ and $R$ tell us exactly when this is applied. We now make this precise.

Start with $S = (r)$, $\widehat{F} = F$, and $j = 0$. If $r \in F$, then set $w_0 = r$, $j = 1$, and remove $r$ from $\widehat{F}$. Throughout the algorithm, we again use $t$ do denote the ordinal number of the current iteration. As long as $S$ is not empty, repeat the following:
\begin{enumerate}[(i)]
    \item Let $s$ be the last vertex in $S$.
    \item If $t \in B$:
    \begin{itemize}
        \item If $j = 0$ or $f_j = 0$, then take $w_j \in \hat F$ to be the closest vertex to $s$, with tie-breaking according to $\pi$, and remove it from $\hat F$. 
        \item  Otherwise, take $w_j$ to be the vertex given by following $\mathbf{c}_j$ from $s$.
    \end{itemize}
    Add $w_j$ to the end of $S$, and increase $j$.
    \item Otherwise, we have $t \notin B$. Remove $s$ from $S$, and if $S$ is empty terminate the procedure. If $t \notin R$, proceed to the next round. Else:
    \begin{enumerate}[(a)]
        \item If $f_j = 0$ (note that we cannot have $j = 0$ at this point), proceed the same as in the original algorithm: Take $u \in T_{s'} \cap \hat F$ for which $h(u)$ is the closest to $s$, and if there are multiple such vertices take the one for which the path from $h(u)$ to $u$ is shortest, tie-breaking according to $\pi$. Let $h_j = h(u)$. 
        \item If $f_j = 1$, then let $h_j$ be the vertex $b_j$ steps back from $s$ towards the root in the tree $T$, and then let $w_j$ be obtained by following $\mathbf{c}_j$ from $h_j$.
    \end{enumerate}
    If $h_j \neq s'$ add $h_j$ to the end of $S$. Add $w_j$ to the end of $S$, remove it from $\hat F$ (relevant only if $w_j$ was obtained in (a)), and increase $j$.
\end{enumerate}
By comparing the two algorithms, we see that they output the same ordering $\sigma(W)$. Therefore, we can uniquely reconstruct $W$ from $\psi_{f,F}(W)$. 
\end{proof}

With Lemma \ref{lemma:key} at hand, the proof of Lemma \ref{lemma:random_branching} is identical to the proof of \cite[Theorem 1.1]{rao17sharp}. We repeat the argument for convenience of the reader.

\begin{proof}[Proof of Lemma \ref{lemma:random_branching}]
    Consider some $x \in V(G)$ and a subset $U \subseteq V(T)$. Given a random $T$-walk $\phi$ on $G$, let $X$ denote the number of vertices $v \in U$ such that $\phi(v) = x$. Our aim is to show
    \begin{equation} \label{eq:moment_gen}
        \mathbb{E}\left[ e^{X} \right] = \sum_{q = 0}^{\infty} \frac{\mathbb{E}[X^q]}{q!} \le e^{K_0 \mathbb{E}[X]},
    \end{equation}
    from which we derive the desired tail bound using Markov's inequality,
    $$
        \Pr\left[ X > K \mathbb{E}[X] \right] = \Pr\left[ e^{X} > e^{K \mathbb{E}[X]} \right] < \mathbb{E}[e^X] e^{-K \mathbb{E}[X]} \stackrel{\eqref{eq:moment_gen}}{\le} e^{-(K - K_0)\mathbb{E}[X]}.
    $$

    Recall the notation of Lemma \ref{lemma:key}: given $W \subseteq V(T)$, let $I_x(W)$ be the indicator random variable for the event that all the vertices in $W$ are mapped to $x$. When $W = \{w\}$, we simply write $I_x(w)$. Note that $X = \sum_{w \in U} I_x(w)$. Consider some $q \in \mathbb{N}$. Then
    $$
        X^q = \left( \sum_{w \in U} I_x(w) \right)^q = \sum_{k = 1}^q \bracenom{q}{k} k! \sum_{W \in \binom{U}{k}} Z_x(W),
    $$
    where $\bracenom{}{}$ denotes the Stirling number of the second kind. By the linearity of expectation, we have
    $$
        \mathbb{E}[X^q] = \sum_{k = 1}^q \bracenom{q}{k} k! \; \mathbb{E}\Bigl[ \sum_{W \in \binom{U}{k}} I_x(W) \Bigr].
    $$
    Combined with Lemma \ref{lemma:key}, this gives the following upper bound on $\mathbb{E}\left[e^X\right]$:
    $$
        \sum_{q = 0}^\infty \frac{\mathbb{E}\left[ X^q \right]}{q!} \le 1 + \sum_{q = 1}^\infty \frac{1}{q!} \sum_{k = 1}^q  \bracenom{q}{k}  k! \sum_{i = 0}^{k-1}\binom{k-1}{i} \frac{(2^8 \mathbb{E}[X])^{k-i}}{(k-i)!} \left( 2^9D \; \frac{\lambda}{d} \right)^i.
    $$
    Rearranging the sums, we get
    \begin{equation} \label{eq:simplify}
        1 + \sum_{i = 1}^\infty \frac{(2^8 \mathbb{E}[X])^{i}}{i!} \sum_{k = i}^\infty \binom{k-1}{i-1} \left( 2^9 D \; \frac{\lambda}{d} \right)^{k-i} \sum_{q = k}^\infty \bracenom{q}{k} \frac{k!}{q!}.
    \end{equation}
    Using the following identity \cite[Eq. 1.94(b)]{Stanley_2011},
    $$
        \sum_{q = k}^\infty \bracenom{q}{k} \frac{1}{q!} = \frac{(e-1)^k}{k!} < \frac{2^k}{k!},
    $$
    we further upper bound \eqref{eq:simplify} as
    $$
        1 + \sum_{i = 1}^\infty \frac{(2^8 \mathbb{E}[X])^i}{i!} 2^i \sum_{k = i}^\infty \binom{k-1}{i-1} \left( 2^9D \; \frac{\lambda}{d} \right)^{k-i}.
    $$
    Now using the identity
    $$
        \sum_{k=i}^\infty \binom{k-1}{i-1} x^{k-i} = \frac{1}{(1 - x)^i}
    $$
    for $|x|<1$, the inner sum further evaluates to 
    $$
      \left(1 - 2^9 D \; \frac{\lambda}{d}\right)^{-i} < 2^i.  
    $$
    We finally get
    $$
        \mathbb{E}\left[e^X\right] < 1 + \sum_{i=1}^\infty \frac{(2^{10} \mathbb{E}[X])^i}{i!} =  e^{2^{10} \mathbb{E}[X]}.
    $$
\end{proof}

\section{Universal hypergraphs}

In this section we prove Theorem \ref{thm:bounded}. In fact, we prove a more general result, Theorem \ref{thm:bounded_density}, on universality of bounded-degree hypergraphs with an additional bound on the density. Universality for such a family of graphs was recently studied by Alon et al.~\cite{alon24density}. Theorem \ref{thm:bounded_density} also improves the bound in \cite[Theorem 1.4]{alon24density}.

Let $\mathcal{H}^{(r)}(q, D, n)$ the family of all $n$-vertex $r$-graphs $H$ with maximum degree at most $D$ and density $m(H) \le q$, where $m(H)$ is as defined in Section \ref{sec:decomposition}. An $r$-uniform hypergraph $H$ with maximum degree $D$ and $v$ vertices contains at most $Dv / r$ edges, thus $m(H) \le D / r$. Therefore $\mathcal{H}^{(r)}(D, n) \subseteq \mathcal{H}^{(r)}(D/r, D, n)$, thus Theorem \ref{thm:bounded_density} implies Theorem \ref{thm:bounded}.

\begin{theorem} \label{thm:bounded_density}
For every $r, D, n \in \mathbb{N}$ and $q \in \mathbb{Q}$, $q > 1 / (r-1)$, there exists $C = C(r, D, q) > 0$ and an $r$-graph $\Gamma$ with
$$
    e(\Gamma) \le C n^{r - 1/q} \log^{1/q}(n)
$$
edges which is $\mathcal{H}^{(r)}(q, D, n)$-universal.
\end{theorem}
\begin{proof}
Fix smallest $a, b \in \mathbb{N}$ such that $b > r-1$ and $q = a/b$. We first describe the construction of $\Gamma$, and then prove that it contains every hypergraph $H \in \cH^{(r)}(q, D, n)$. We say that $H$ is an $(\le r)$-uniform hypergraph, or $(\le r)$-graph for short, if every hyperedge in $H$ has size at most $r$. For brevity, throughout the proof we ignore ceilings and floors.

\paragraph{Construction.} Let $m = (n / \log (n))^{1/a}$, and let $d \in \mathbb{N}$ be sufficiently large with respect to $D$. Let $G$ be an $(m, d, \lambda)$-graph on the vertex set $[m]$, for some $\lambda < \alpha d / D$ where $\alpha$ is as given by Lemma \ref{lemma:random_branching}. Explicit construction of such a graph, for any $m \ge m_0(d)$, was obtained by Alon \cite{alon21expander}. Let $G^2$ be the graph obtained from $G$ by adding an edge between every two vertices at distance at most $2$ in $G$. Note that the maximum degree in $G^2$ is at most $d^2$.

We first form an $(\le r)$-graph $\Gamma'$ as follows: $V(\Gamma') = [m]^a$, and $r' \le r$ vertices $\mathbf{v}^{(1)}, \ldots, \mathbf{v}^{(r')} \in [m]^a$, $\mathbf{v}^{(j)} = (v_1^{(j)}, \ldots, v_a^{(j)})$, form a hyperedge if there exist forests $F_1, \ldots, F_a$, each on the vertex set $[r]$, such that:
\begin{itemize}
    \item $\sum_{i=1}^a e(F_i) \ge b$, and
    \item for each $i \in [a]$ there exists a homomorphism $f_i \colon [r] \to \{v_1^{(i)}, \ldots, v_{r'}^{(i)}\}$ of $F_i$ in $G^2$.
\end{itemize}
This construction is guided by the statement of Lemma \ref{lemma:hyp_decomposition}.

Next, take $\Gamma$ to be an $r$-graph obtained as the $(C \log n)$-blowup of $\Gamma'$, for $C$ being a sufficiently large constant. That is, for each vertex $\mathbf{v} \in V(\Gamma')$ we introduce a set $B_{\mathbf{v}}$ of size $C \log n$, and for each hyperedge $\mathbf{v}^{(1)} \cdots \mathbf{v}^{(r')} \in \Gamma'$, for some $1 \le r' \le r$, we add all the subsets of $\bigcup_{i = 1}^{r'} B_{\mathbf{v}^{(i)}}$ of size exactly $r$ as hyperedges in $\Gamma$.

Let us first count the number of edges in $\Gamma'$. Consider forests $F_1, \ldots, F_a$ on $[r]$, such that $\sum_{i=1}^a e(F_i) \ge b$. Then $c(F_1) + \ldots + c(F_a) \le r a - b$, where $c(F)$ denotes the number of connected components in $F$. As the maximum degree in $G^2$ is $d^2$, a homomorphism of a forest $F$ in $G^2$ can be chosen in at most $m^{c(F)} (d^2)^{r - c(F)}$ ways. Altogether, this gives at most
$$
    \sum_{F_1, \ldots, F_a} \prod_{i = 1}^{a} m^{c(F_i)} (d^2)^{r - c(F_i)} \le r^{r a} m^{ra - b} (d^2)^{ra} = r^a d^{2 ra} \left( \frac{n}{\log n} \right)^{r - b/a} 
$$
hyperedges in $\Gamma'$, where the first sum goes over $a$-tuples $(F_1, \ldots, F_a)$ of forests on $[r]$ such that $c(F_1) + \ldots + c(F_a) \le r$.

Each hyperedge in $\Gamma'$ gives rise to less than $(r C \log n)^r$ hyperedges in $\Gamma$. Therefore, $\Gamma$ has
$$
    O\left( n^{r - 1/q} \log^{1/q}(n) \right)
$$
hyperedges.

\paragraph{Universality.} We now show that $\Gamma$ is $\cH(q, D, n)$-universal. Consider an $r$-graph $H \in \cH^{(r)}(d, D, n)$. Let $H_1, \ldots, H_a$ be graphs on the vertex set $V(H)$ given by Lemma \ref{lemma:hyp_decomposition}, and for each edge $h \in H$ let $F_1^{(h)} \subseteq H_1, \ldots, F_a^{(h)} \subseteq H_a$ be forests corresponding to \ref{prop:forest}. For each $i \in [a]$, let $T_i$ be a forest on $V(H)$ with maximum degree at most $2D$ such that $H_i \subseteq T_i^2$, that is, if $\{v,w\}$ is an edge in $H_i$ then $\{v,w\}$ are at distance at most $2$ in $T_i$. Such $T_i$ can be obtained from $H_i$ by replacing each cycle $x_1, \ldots, x_\ell$ in $H_i$ by a path $x_1, x_\ell, x_2, x_{\ell-1}, x_3, x_{\ell-2}, \ldots$. We use the fact that if $f \colon T_i \to G$ is a homomorphism of $T_i$ in $G$, then it is also a homomorphism of $H_i$ in $G^2$. Therefore, from now on we can focus on $T_i$ instead of $H_i$. By adding more edges, without loss of generality we may assume each $T_i$ is a tree.

Our aim is to iteratively find homomorphisms $\phi_i \colon T_i \to G$ such that, for each $\mathbf{v} \in [m]^i$, the set
$$
    S_{\mathbf{v}}^i = \left\{ w \in V(H) \colon \phi_1(w) = v_1, \ldots, \phi_i(w) = v_i \right\}
$$
is of size 
\begin{equation} \label{eq:size_S}
    |S_\mathbf{v}^i| \le n (K / m)^i,
\end{equation}
where $K \ge K_0$ is sufficiently large and $K_0$ is the constant given by Lemma \ref{lemma:random_branching}. Suppose that we have found such homomorphisms $\phi_1, \ldots, \phi_{i-1}$, for some $i \in \{1, \ldots, a\}$. For simplicity, we define $S_{\mathbf{v}}^0 = V(H)$. Let $\phi_i$ be a random $T_i$-walk in $G$, with an arbitrary vertex in $T_i$ being the root. Applying Lemma \ref{lemma:random_branching} with some $x \in V(G)$ and $S = S_\mathbf{v}^{i-1} \cup S'$, for some $\mathbf{v} \in [m]^{i-1}$ and $S'$ chosen arbitrarily such that $|S| = n(K / m)^{i-1}$ (this is a rather technical detail), we get
$$
    \Pr[ X(x,\mathbf{v}) > K |S| / m ] \le e^{-(K - K_0) |S| / m} = o(1/n^2),
$$
where $X(x,\mathbf{v})$ denotes the number of vertices $v \in S$ such that $\phi_i(v) = x$. There are $m$ choices for $x$ and $m^a = o(n)$ choices for $\mathbf{v}$, thus with positive probability $\phi_i$ is such that $X(x, \mathbf{v}) \le n(K / m)^i$ for every $x \in V(G)$ and $\mathbf{v} \in [m]^{i-1}$. Therefore, there exists $\phi_i$ for which \eqref{eq:size_S} holds.

Let now $f \colon H \to [m]^a$ be defined as $f(w) = (\phi_1(w), \ldots, \phi_a(w))$. As $\phi_i$ is a homomorphism of $T_i$ into $G$ and $H_i \subseteq T_i^2$, $\phi_i$ is also a homomorphism of $H_i$ into $G^2$. Consider some edge $h = \{w_1, \ldots, w_r\} \in H$. Then the restriction of $\phi_i$ to $F_i^{(h)}$ is a homomorphism of $F_i^{(h)}$ in $G^2$. Since $\sum_{i=1}^a e(F_i) \ge b$ by Lemma \ref{lemma:hyp_decomposition}, the set of vertices $\mathbf{v}^{(1)}, \ldots, \mathbf{v}^{(r)}$ given by $\mathbf{v}^{(i)} = (\phi_1(w_i), \ldots,\phi_a(w_i))$ for $i \in [r]$, form a hyperedge in $\Gamma'$. Note that these vertices might not all be distinct, thus the hyperedge is not necessarily of size $r$. Therefore, each $h \in H$ is preserved by $f$, that is, $f(h)$ is a hyperedge in $\Gamma'$. By \eqref{eq:size_S}, no vertex is an image of more than $O(\log n)$ vertices from $H$, thus we can turn $f$ into an injective homomorphism $f' \colon H \to \Gamma$, which finally defines a copy of $H$ in $\Gamma$.
\end{proof}

A reader familiar with some of the previous work \cite{alon07sparsebounded,alon08optimal,alon24density}, most notably the proof of \cite[Theorem 1.4]{alon24density}, will notice similarity with the overall strategy used in the proof of Theorem \ref{thm:bounded_density}. We note that both the construction of the hypergraph $\Gamma$ and tools used to show it is universal are more involved in our case.

\section{Concluding remarks}

\paragraph{Towards optimal universal hypergraphs.} The following conjecture, if true, would allow one to replace the use of Lemma \ref{lemma:random_branching} in the proof of Theorem \ref{thm:bounded} and obtain optimal the optimal bound $O(n^{r - r/d})$. 

\begin{conjecture}
    For every $D \in \mathbb{N}$ there exists $d \in \mathbb{N}$ and $\alpha > 0$ such that the following holds. Suppose $G$ is an $(n, d, \lambda)$-graph for $\lambda = \Theta(\sqrt{d})$. Then for any tree $T$ with maximum degree at most $D$ and any partition of $V(T) = T_1 \cup \ldots \cup T_\ell$ with $|T_i| \le \alpha n$ for each $I \in [\ell]$, there exists a homomorphism $\phi \colon T \to G$ such that the restriction of $\phi$ to each $T_i$ is an injection.
\end{conjecture}

In case $\ell = 1$, this corresponds to the celebrated result by Friedman and Pippenger \cite{friedman87tree}. It is important to notice that there is no upper bound on $\ell$ in terms of $n$ or any other parameter.

\paragraph{Branching random walks.} As discussed in Section \ref{sec:branching_walk_intro}, it would be interesting to extend Lemma \ref{lemma:random_branching} to other forms of random variables. As a first step, one could consider a random variable which counts the number of times a branching random walk has landed in a set $X \subseteq V(G)$ rather than on a particular vertex $x \in V(G)$. This would generalise the result of Gillman \cite{gilam98chernoff}. 

The main obstacle in adapting the proof of Lemma \ref{lemma:random_branching} to this case seems to be in the estimate \eqref{eq:expectation_W}. For appropriately defined random variable, analogous to the one considered in \eqref{eq:expectation_W}, one would aim to upper bound its expectation as $\mu \prod_{i=1}^{k-1} (\mu + (\lambda/d)^{d_i})$, where $\mu = |X| / n$. When $|X| = 1$ then this is precisely what is achieved in \eqref{eq:expectation_W}. However, the way this upper bound is derived is based on the `worst-case' conditioning of how a part of $T$ is mapped, which is too pessimistic and does not yield a desired upper bound for larger set $X$. A similar estimate is encountered \cite[Lemma 3.3]{rao17sharp}, however the proof there leverages the linear structure of random walks and the ability to algebraically express the desired expectation, and it is not clear whether a similar argument can be applied here.

Another interesting direction would be to consider branching random walks on graphs which are not regular.

\paragraph{Acknowledgment.} The author thanks Noga Alon for stimulating discussions about the problem, and Anders Martinsson for discussions about Lemma \ref{lemma:random_branching}.

\bibliographystyle{abbrv}
\bibliography{universality}

\end{document}